\documentstyle [twoside,12pt,amsmath,amssymb,amsthm]{article}
\newtheorem{thm}{Theorem}[section]
\newtheorem{cor}[thm]{Corollary}
\newtheorem{lem}[thm]{Lemma}
\newtheorem{prop}[thm]{Proposition}

\newtheorem{remark}[thm]{Remark}

\newcommand{\thmref}[1]{Theorem~\ref{#1}}
\newcommand{\lemref}[1]{Lemma~\ref{#1}}
\newcommand{\propref}[1]{Proposition~\ref{#1}}

\newcommand{\secref}[1]{Section~\ref{#1}}

\newcommand{\bx}{\hfill$\Box$\vspace{.6cm}}

\numberwithin{equation}{section}

\renewcommand\a{\alpha}         
\renewcommand\b{\beta}
\newcommand\g{\gamma}
\renewcommand\d{\delta}
\newcommand\e{\epsilon}
\renewcommand\l{\lambda}
\renewcommand\L{\Lambda}
\newcommand\D{\Delta}
\newcommand\G{\Gamma}
\newcommand\f{\frac}
\newcommand{\Z}{{\mathbb{Z}}}
\newcommand{\R}{{\mathbb{R}}}
\newcommand{\C}{{\mathbb{C}}}

\newcommand{\U}{{\mathbb{H}}}
\newcommand\re{\mbox{Re~}}
\newcommand\im{\mbox{Im~}}
\renewcommand\Re{\mbox{Re~}}
\renewcommand\Im{\mbox{Im~}}

\newcommand\F{{\cal F}}
\renewcommand\H{{\cal H}}

\newcommand{\ttwo}[4]{\left(\begin{array}{cc}
{#1} & {#2} \\ {#3} & {#4} \end{array} \right)}
\newcommand{\ttt}[9]{\left(\begin{array}{ccc} {#1} & {#2} & {#3} \\
    {#4} & {#5} & {#6} \\   {#7} & {#8} & {#9} \end{array}  \right)}

\renewcommand\i{^{-1}}

\renewcommand\({\left(}         
\renewcommand\){\right)}

\begin{document}

\title{On the Existence and Temperedness of Cusp Forms
for $SL_3({\Z})$}
\author{Stephen D. Miller\thanks{The author was supported by
 an NSF
Postdoctoral Fellowship during this work.} \\ Department of
Mathematics \\ Yale University \\ P.O. Box 208283 \\ New Haven, CT
06520-8283 \\{\tt stephen.miller@yale.edu}}
\date{}
 \maketitle

\begin{abstract}
We develop a partial trace formula which circumvents some
technical difficulties in computing the Selberg trace formula for
the quotient $SL_3({\Z})\backslash SL_3({\R})/SO_3({\R})$.  As
applications, we establish the Weyl asymptotic law for the
discrete Laplace spectrum and prove that almost all of its cusp
forms are tempered at infinity.  The technique shows there are
non-lifted cusp forms on $SL_3({\Z})\backslash
SL_3({\R})/SO_3({\R})$ as well as non-self-dual ones.  A
self-contained description of our proof for $SL_2({\Z})\backslash
\U$ is included to convey the main new ideas.  Heavy use is made
of truncation and the Maass-Selberg relations.
\end{abstract}

\section{Introduction}

    In the 1950s A. Selberg (\cite{Selberg})
developed his trace formula to prove the existence of
non-holomorphic, everywhere-unramified, cuspidal ``Maass" forms.
These are real-valued functions on the upper-half plane ${\U}=
\{x+iy\mid y>0\}$ which are invariant under the
 action of $SL_2({\Z})$ by fractional linear transformations.
Unlike the holomorphic cusp forms, which can all be explicitly
described, no Maass form for $SL_2({\Z})$ has ever been
constructed and they are believed to be intrinsically
transcendental.

The non-constant Laplace eigenfunctions in
$L^2(SL_2({\Z})\backslash \U)$ are all Maass forms.  Since
$SL_2({\Z})\backslash \U$ is noncompact, their existence is no
triviality; they very likely do not exist on the generic
finite-volume quotient of $\U$ (see \cite{Sarnak}). However,
Selberg showed that for any congruence subgroup $\G \subset
SL_2({\Z})$, the discrete spectrum is as large as one can expect
-- namely, it obeys the same asymptotics as the spectrum of a
compact surface of the same size:

\begin{thm}(Selberg)
\label{selweyl} Let $0=\l_0<\l_1\le \l_2 \le \cdots$ be the
discrete eigenvalues (with multiplicity) of the non-Euclidean
laplacian $\D$ on $\G\backslash\U$. If $\G\subset SL_2({\Z})$ is a
congruence subgroup, then \begin{equation} N(T)=\#\{\l_j \le T \}
\sim \f{ \mbox{area}(\G\backslash\U)}{4\pi}T \end{equation}
 as
$T\rightarrow\infty$.
\end{thm}

Similar asymptotics hold for any compact manifold by theorems of
H. Weyl and others; we will refer to such an asymptotic for $N(T)$
as the {\em Weyl law} for the space in question.  For noncompact
manifolds one can count the discrete and continuous spectra
together, but it is extremely difficult to decouple the two. The
main contribution of this paper is a technical novelty for
separating them and proving the Weyl law for certain noncompact
quotients.

One can generalize Maass forms to groups other than $SL_2({\R})$
 in a variety of ways.  We are most interested in the linear group
 and so will focus our attention there,
to $G=SL_3({\R})$ in particular.  Let $K=SO_3({\R})$ be a maximal
compact subgroup of $G$, $\H=G/K, \G=SL_3({\Z}),$ and
$X=\G\backslash\H$. The ring $\cal D$ of ($G$-)invariant
differential operators on the symmetric space $\H$ is a polynomial
ring in two generators; it is explicitly described in \cite{Bump},
p. 32 for example.  We will concentrate first on a particular
element of $\cal D$, the laplacian $\D$, which is normalized so
that its continuous spectrum on $L^2(\H)$ is the interval
$[1,\infty)$. Here again, the non-constant discrete eigenfunctions
of $\D$ in $L^2(X)$ are cusp forms (see \secref{struct} for
definitions), which are the appropriate generalizations of Maass
forms and our basic objects of study.

By symmetry considerations, it is not difficult to prove the
existence of {\it odd} cusp forms on $SL_2({\Z})\backslash {\U}$
and $SL_3({\Z})\backslash \H$. These comprise only half of the
expected spectrum; the deeper issue is the existence of {\it even}
cusp forms. The only other known cusp forms on
$SL_3({\Z})\backslash \H$ are Gelbart-Jacquet lifts
(\cite{GelJac}) of the forms on $SL_2({\Z})\backslash {\U}$
 that Selberg
discovered (many of which have been numerically identified--see
\cite{Hejhal}). However, no version of the trace formula has been
used for the higher rank $SL_3({\Z})\backslash \H$ to prove the
Weyl law for cusp forms, {\it ala} Selberg's
\thmref{selweyl}.\footnote{In \cite{Stade} a proof of
\thmref{mainresult} is claimed.
  This proof appears to be incorrect because it is based
on taking a specifically-chosen sequence of functions in
 the trace formula of \cite{Wallace}.  Unfortunately that
trace formula is incomplete, because in its estimation of the
continuous spectrum, it ignores the poles of the intertwining
 operators which arise from contour shifts.  This seems to be a complicated
problem to fix, though doing so would lead to much more specific
information about the discrete spectrum (e.g. perhaps an error
term).}  For example, Arthur's trace formula (\cite{ArthurICM})
computes the traces of certain integral operators over the
discrete spectrum, but in it the discrete spectrum is paired with
parabolic orbital integrals. Though it appears difficult, it would
be very interesting to separate the two and give estimates on the
size of the spectrum -- especially because Arthur's formula has
been developed for general quotients.

We will shortcut the trace formula to prove the Weyl law for
$SL_3({\Z})\backslash \H$ as well as some qualitative results
about the spectrum.

\begin{thm}
\label{mainresult} Let $0=\l_0<\l_1\le\l_2\le\cdots$ denote the
eigenvalues of $\D$ on $SL_3({\Z})\backslash\H$ with multiplicity.
Then
\begin{equation}
 N(T):=\#\{\l_j\le T\}\sim \f{vol(SL_3({\Z})\backslash \H
)}{\G(7/2)} \(\f{T}{4 \pi}\)^{5/2}
\end{equation}
 as $T\rightarrow \infty$.
\end{thm}  These are the same asymptotics that the spectra of
closed five-dimensional manifolds obey.  This result confirms a
conjecture of Sarnak (\cite{Sarnak}),
 who asserted
that the Weyl law holds for the cuspidal spectrum of the laplacian
on any congruence quotient of $SL_n({\R})/SO_n({\R})$.
\thmref{mainresult} is actually a special case of the more-general
\thmref{genct}, a kind of equidistribution theorem which counts
the two-dimensional joint spectrum of $\cal D$ lying in various
sets. This generalization is used to establish many of the other
findings below, so first we shall briefly describe the joint
spectrum before stating the applications.

The ring of invariant differential operators $\cal D$ on $\H$ is
commutative and we may take an orthonormal set of common
eigenfunctions
\begin{equation}
\f{1}{\sqrt{\mbox{area}(SL_3({\Z})\backslash
\H)}}=\phi_0,\phi_1,\phi_2,\ldots \end{equation} such that
$$\D\phi_j=\l_j\phi_j.$$  Each $\phi_j$ in turn induces a
homomorphism

$$\l:{\cal D}\rightarrow
{\C},~~\mbox{where}~D\phi=\l(D)\phi,~D\in{\cal D}.$$  This map is
best described by ``principal series'' or ``Langlands'' parameters
$$\{\ell_1,\ell_2,\ell_3\},~~~ \ell_1+\ell_2+\ell_3=0$$ (see
\secref{convopsec} for more details).  It is possible to compute
$\l(D)$ in terms of these spectral parameters, e.g.
\begin{equation}
\l(\D)=1-\f{\ell_1^2+\ell_2^2+\ell_3^2}{2}.\end{equation}  The
real parts of the $\ell_i$ are bounded, for example by a trivial
bound from representation theory:
\begin{equation} |\Re{\ell_i}| < \f12~~~(\cite{JacSha}).
\label{introjs}\end{equation} Hence the Weyl law essentially
counts the number of spectral parameters $(\ell_1,\ell_2,\ell_3)$
lying in a large ball.  In \thmref{genct} we compute the
asymptotics of how many points of the joint spectrum lie in
certain sets of various shapes and sizes. Along with some
representation theory, we can count how many cusp forms obey given
properties. The statement of \thmref{genct} is technical so we
will only describe its corollaries here.

A cusp form is called {\em tempered} if each $\Re\ell_i=0$.  The
archimedean Ramanujan-Selberg conjecture asserts that all cusp
forms on $SL_3({\Z})\backslash \H$ are tempered.  The only
nontrivial improvement of (\ref{introjs}) on individual forms is
the result of \cite{LRS}.

\begin{thm}
\label{aeram}
  Almost all of the cusp forms
on $SL_3({\Z})\backslash SL_3({\R})/SO_3({\R})$ obey
 the archimedean Ramanujan-Selberg conjectures in the sense
that
$$\lim_{T\rightarrow\infty}
\f{\#\{\l_j \le T \mid \D\phi_j=
\l_j\phi_j, \phi_j \mbox{ tempered}\}}{\#\{\l_j \le T\}} = 1.$$
\end{thm}

For $SL_2({\Z})\backslash \U$ the nontempered cusp forms are those
with Laplace eigenvalue less than $1/4$, and the Weyl law shows
there are only a finite number of these.  In higher rank, though,
there may be nontempered forms with large Laplace eigenvalues. The
archimedean Ramanujan-Selberg conjecture is an analog of Selberg's
famous ``1/4 conjecture,'' as it automatically implies that every
cuspidal Laplace eigenvalue satisfies
$$\l^{cusp}(SL_3({\Z})\backslash \H) \ge \l_1(\H)=1.$$ In fact,
this has already been proved in \cite{Miller}, but without
establishing temperedness itself. The technique here is different
and further-reaching in that it works whenever the Weyl law can be
proved by our method.

\begin{cor}
If $$N^{cusp}(T):=\#\{\l_j \le T \mid \D\phi_j = \l_j\phi_j \mbox{
for a nonzero cusp form } \phi_j \},$$ then $$N^{cusp}(T)\sim
N(T)\sim \f{vol(SL_3({\Z})\backslash \H )}{\G(7/2)} \(\f{T}{4
\pi}\)^{5/2}$$ as $T\rightarrow \infty$.
\end{cor}
{\bf Proof}: The only other eigenfunctions are residues of Eisenstein series,
which are not tempered, and \thmref{aeram} shows these have measure zero.

Another and more-precise reason uses the classification of the
discrete spectrum established by M\oe glin and Waldspurger
(\cite{MoeWal}, confirming the conjecture of \cite{Jacquet}). A
consequence is that the only non-cuspidal discrete eigenfunctions
of $\D$ on $SL_p({\Z})\backslash SL_p({\R})/SO_p({\R})$ are
constants when $p$ is prime. Thus for $T\ge 0$,
$N(T)=N^{cusp}(T)+1$, as alluded to earlier.\bx

\begin{cor}
\label{native}
There are ``native" cusp forms on
 $SL_3({\Z})\backslash\H$
which are not Gelbart-Jacquet lifts of Maass forms from
$SL_2({\Z})\backslash {\U}$. Moreover, these native forms comprise
100\% of the spectrum in the sense that
$$\lim_{T\rightarrow\infty} \f{\#\{\l_j \le T \mid \D\phi_j=
\l_j\phi_j, \phi_j \mbox{ a Gelbart-Jacquet lift}\}}{\#\{\l_j \le
T\}} = 0.$$
\end{cor}

A cusp form is called {\em self-dual} if $\phi(g)=\phi((g^t)\i)$
for all $g\in G$.  Alternatively, its joint spectral parameters
are $$\{\ell_1,\ell_2,\ell_3\}=\{\mu,0,-\mu\}\mbox{ for some
}\mu\in{\C}.$$
\begin{thm}
\label{nonselfdual}
 There exists a non-self-dual cusp form on
$SL_3({\Z})\backslash SL_3({\R})/SO_3({\R})$, and in fact these
are also of full measure: $$\lim_{T\rightarrow\infty} \f{\#\{\l_j
\le T \mid \D\phi_j= \l_j\phi_j, \phi_j \mbox{
self-dual}\}}{\#\{\l_j \le T\}} = 0.$$
\end{thm}
All Gelbart-Jacquet lifts are self-dual and it is conjectured that
the converse is true.

In the next section we present a self-contained account of the
main ideas of this paper, but specialized towards giving a new
proof of Selberg's \thmref{selweyl} for $\G=SL_2({\Z})$. The main
improvement is a circumvention of the usual derivation of the
trace formula. Usually one integrates the automorphic kernel only
over a truncated fundamental domain, in order to prevent
divergence from the parabolic orbital integrals and the Eisenstein
series. After matching the growth rates of these two terms, a
trace formula is obtained in the limit as the truncation parameter
moves to infinity through the cusp.  Here we merely truncate, and
by positivity considerations arrive at an inequality giving
 a {\em lower} bound on the size of the spectrum.  This gives the Weyl lower
 bound, but with the wrong constant -- it involves the area of the truncated
fundamental domain instead.  By pushing the truncation parameter
to infinity only at the end, we recover the correct lower bound;
the upper bound has already been established by \cite{Donnelly}.

\subsection*{Acknowledgements}
    I wish to thank Peter Sarnak, under whose direction some of this
material first appeared in my Princeton University dissertation.
Also I am indebted to Donna Belli, Don Blasius,
   Sol Friedberg,
Serge Lang, Alex Lubotzky, Jonathan Rogawski, Ze'ev Rudnick,
 Nolan Wallach, Gregg Zuckerman, and the referee for helpful comments.
  This work was supported by NSF graduate and
post-doctoral fellowships, the Yale Hellman fund, and NSA grant
MDA904-99-1-0046.

%
%



\section{The Weyl Law for $SL_2({\Z})\backslash \U$}
\label{appen}

This is a self-contained section illustrating our method on
$SL_2({\Z})\backslash\U$.  We will reprove the Weyl Law
(\thmref{selweyl}) using the same ``partial trace'' technique we
use in higher rank for \thmref{mainresult}.  We stress that the
reason the partial trace is successful has little to do with the
complications of higher rank {\em per se}, but rather because it
serves as a substitute for the trace formula in counting the
spectrum.  As of yet, it appears difficult to give a direct proof
of the Weyl law in higher rank using a trace formula, though this
is how and why Selberg proved it for $SL_2({\Z})\backslash\U$ to
begin with.  Our technique is also applicable to congruence covers
of $SL_2({\Z})\backslash\U$, but we will focus on the full-level
situation here in order to demonstrate how the main difficulty --
non-compactness -- is addressed.  Basic references for this
section are \cite{Selberg}, \cite{Selgott}, \cite{Hejhal}, and
\cite{Terras}.


\subsection{Definitions}

    The hyperbolic plane can be modeled by the complex
upper-half plane $${\U}=\{u+iv\mid v>0\}$$ with area element
$$dA=\f{dudv}{v^2},$$ line element
$$d\sigma^2=\f{du^2+dv^2}{v^2},$$ and laplacian
$$\D=-v^2\(\f{\partial^2}{\partial u^2} + \f{\partial^2}{\partial
v^2}\) .$$

The group $G=SL_2({\R})$ acts by isometries on $\U$ by fractional
linear transformations $$\ttwo abcd :z\mapsto \f{az+b}{cz+d}.$$
The stabilizer of any point in $\U$ is a maximal compact subgroup
isomorphic to $K=SO_2({\R})=\hbox{Stab}(i)$, and so $\U$ is
isomorphic to the symmetric space $G/K$.

Let $\G=SL_2({\Z}), \bar{\G}=PSL_2({\Z})$ and
$X=\G\backslash\U=\bar{\G}\backslash\U$.
 The
space $X$ is noncompact but has finite volume, $\f{\pi}{3}$.  Thus
the laplacian $\D$ has a continuous spectrum on $L^2(X)$ furnished
by Eisenstein series
$$E_s(z)=\f{1}{2}\sum_{\left\{\ttwo{1}{n}{0}{1}\mid n \in {{\Z}}
\right\}\backslash \G}Im(\g z)^s.$$ (This infinite series is only
defined for $Re (s)>1$, where it converges absolutely, but it has
a meromorphic continuation to $s\in\C$.)  The continuous spectrum
comes from $$E_{\f{1}{2}+it}(z)~~,~~t \in {\R}.$$

Also, since $\D$ annihilates constant functions, it has a discrete
spectrum.  Our aim here is to count it.  Let
$$0=\l_0<\l_1\le\l_2\le \cdots$$ be the discrete spectrum of $\D$
on $L^2(X)$ and $$\sqrt{\f{3}{\pi}} =
\phi_0,\phi_1,\phi_2,\ldots$$ be an orthonormal set of
eigenfunctions $$\phi_j\in L^2(X),~~~\D\phi_j=\l_j\phi_j.$$

We will prove
\begin{thm}(=\thmref{selweyl})(Selberg)
As $T\rightarrow\infty$ $$N(T)=\#\{\l_j\le T\}\sim \f{T}{12}.$$
\label{selbergsweylslaw}
\end{thm}

\subsection{Convolution Operators}

Every $g\in C_c^\infty(K\backslash {\U})=C_c^\infty(K\backslash
G/K)$ acts on $f\in L^2(G/K)$ by convolution:

\begin{equation}
(L_gf)(x)= (f*g)(x)= \int_{G/K}f(y)g(y\i x)dy. \label{convongmodk}
\end{equation}
Bi-$K$-invariant functions are functions of the distance of a
point in $\U$ to $i$.  Their convolution operators form a
commutative algebra which also commutes with the Laplace operator.
These operators share common eigenfunctions as well, so Laplace
eigenfunctions on $\U=G/K$ play a special role. The simplest
examples are furnished by power functions

\begin{equation}
\D v^s=s(1-s)v^s \label{powerfunctions}
\end{equation}
and by the bi-$K$-invariant ``spherical'' functions
$$\tilde{\phi}_s(z) = \int_K \im(k z)^sdk,~~~\int_Kdk=1,$$
\begin{equation}
\D \tilde{\phi}_s(z)=s(1-s)\tilde{\phi}_s.
 \label{spherfundef}
\end{equation}
The integral operators $L_g$ also act on $f\in L^2(\G\backslash
G/K)$:
\begin{eqnarray}
(L_g f)(x)= & \int_{G/K}f(y) g(y\i x)dy & \nonumber \\ = &
\sum_{\g\in\bar{\G}}\int_{\G\backslash G/K} f(\g\i y)g(y\i \g x)dy
& \nonumber \\ = & \int_{\G\backslash G/K} f(y) K(x,y) dy, &
\label{convonquotient}
\end{eqnarray}
where
\begin{equation}
K(x,y)=\sum_{\g\in\bar{\G}} g(x\i \g y) \label{autokern}
\end{equation}
is called the {\em automorphic kernel}.  Again, the $L_g$ commute
with each other and with $\D$, so we may assume that the $\phi_j$
are taken to be an orthonormal set of eigenfunctions of all $L_g,
g\in C_c^\infty(K\backslash G/K)$ as well:
\begin{equation}
(L_g\phi)(x) = \int_{G/K}\phi(y)g(y\i x)dy = \hat{g}(\phi)\cdot
\phi(x).
 \label{lgeigfunc}
\end{equation}
Selberg showed that $$\int_K
\phi_j(kg)dk=\phi_j(I)\tilde{\phi}_{s_j}(g),~~~~\D\phi_j=s_j(1-s_j)\phi_j$$
and hence by averaging (\ref{lgeigfunc}) over $K$ he found a
formula for $\hat{g}(\phi)$:

\begin{prop}(Selberg's Uniqueness Principle)
The $L_g$-eigenvalue $\hat{g}(\phi_j)$ of any $\phi_j$ as above
depends only on its Laplace eigenvalue.  Namely, if
\begin{equation}
\D\phi_j=(\f14 - \nu_j^2)\phi_j, \label{onequarminusnus}
\end{equation}
then $$L_g\phi_j = \hat{g}(\nu_j)\phi_j,$$ where the Selberg
transform \begin{equation} \hat{g}(\nu)= \int_\U
g(u+iv)v^{1/2+\nu}\f{dudv}{v^2}. \label{defofselbergtransform}
\end{equation}
\end{prop}

\begin{remark} {\em Formula (\ref{defofselbergtransform}) is in fact
(\ref{lgeigfunc}) applied to the eigenfunction $v^{1/2+\nu}$, and
evaluated at the identity.} \end{remark}

The kernel $K(x,y)$ thus has two expansions:

\begin{prop}(Spectral Expansion) For $g\in C_c^\infty(K\backslash G/K)$,
$$K(x,y)=\sum_{\g\in\bar{\G}} g(y\i \g x) = $$
\begin{equation}
\sum_{\stackrel{j\ge 0}{\mbox{disc. $L^2$
spec.}}}\hat{g}(\nu_j)\phi_j(x)\phi_j(y) + \f{1}{4\pi}\int_{{\R}}
\hat{g}(it) E_{1/2+it}(x)\overline{E_{1/2+it}(y)}dt.
 \label{specexpans}
\end{equation}
\end{prop}

We will now recall some analytic properties of the Selberg
transform that will be needed in constructing our choices of $g$
later.  First is the inversion formula \begin{equation}
g(x)=\int_{{\R}} \hat{g}(it)\tilde{\phi}_{1/2-it}(x)\f{t\tanh{\pi
t}}{4\pi} dt. \label{spherinver}
\end{equation}
Secondly, we need to clarify the relationship between the Selberg
transform and the Mellin transform $\cal M$.  If
\begin{equation}
\bar{g}(v)=\int_{{\R}} g(u+iv)du \label{harishtransf}
\end{equation}
denotes the ``Harish transform,'' then
(\ref{defofselbergtransform}) expresses
\begin{equation}\hat{g}(\nu) = {\cal
M}\(\f{\bar{g}(x)}{\sqrt{x}}\)(\nu).\label{mellsel}\end{equation}
Through the Mellin inversion formula, formula (\ref{spherinver})
essentially amounts to a statement about the Harish transform.  A
more-precise fact was proven by Ehrenpreis and Mautner (see also
\cite{Gangolli}).

\begin{thm}(\cite{EhrMau}).
\label{ehren}  The Harish transform $g\mapsto\bar{g}$ is a
bijection between

$$C_\sigma^\infty(K\backslash\U) = \{ g \in
C_c^{\infty}(K\backslash \U) \mid supp(g) \subset K\cdot
i(\sigma\i,\sigma) \}$$

and $$C_\sigma^\infty(0,\infty) = \{h \in C_c^\infty(0,\infty)
\mid h(v)=h(v\i), supp(h)\subset(\sigma\i,\sigma) \}.$$
\end{thm}

\subsection{Construction of Functions}

Recall that there exist smooth, non-negative functions on $\R$
with arbitrarily-small support whose Fourier transforms are also
non-negative.  (The easiest way to make one is to convolve a
non-negative, compactly supported function with itself and
rescale.)  Using (\ref{mellsel}) and \thmref{ehren}, we may fix
some $g\in C_\sigma^\infty(K\backslash \U)$ such that $g\ge 0$,
$\int_{{\R}} \hat{g}(it)dt=1$, and $\hat{g}(it)\ge 0$ for
$t\in\R$. We would like to scale $\hat{g}$ so that it resembles
the characteristic function of a large interval.  However, in this
non-euclidean setting it is difficult to control $g$ and $\hat{g}$
simultaneously.  Instead, we will convolve $\hat{g}(it)$ with
$\chi_{[-T,T]}$:
\begin{equation}
\hat{g}_T(it)=\int_{-T}^T\hat{g}(it+ir)dr ~\ge 0, \label{ghattdef}
\end{equation}
which is the Mellin transform of
$$\f{\bar{g}(x)}{\sqrt{x}}\int_{-T}^Tx^{ir}dr ~\in
C_\sigma^\infty(0,\infty).$$  Hence \thmref{ehren} also guarantees
the existence of some function $g_T \in
C_\sigma^\infty(K\backslash \U)$ whose Selberg transform
$$\widehat{g_T}=\hat{g}_T.$$

  We will use the $g_T$'s in the
automorphic kernel.  Our construction allows us to conclude two
important analytic properties:

\begin{prop}
\label{maxgt} $$\max_{x\in\U} |g_T(x)| = g_T(i).$$
\end{prop}
{\bf Proof:} This follows from the inversion formula
(\ref{spherinver}), the positivity of $\hat{g}_T$
(\ref{ghattdef}), and the inequality $$|\tilde{\phi}_{1/2+it}(x)|
\le 1 = \tilde{\phi}_{1/2+it}(i)$$ (see \cite{Kolk}, (2.13)).  \bx

\begin{prop}
\label{gtandchar} For any $m\ge 0$, $$\hat{g}_T(it) =
  \begin{cases}
     1+O_m\((|T|-|t|+1)^{-m}\) & t < T, \\
    O_m\((|t|-|T|+1)^{-m}\) & t \ge T.
  \end{cases}
$$
\end{prop}
{\bf Proof:} Since $g$ is smooth, $\hat{g}(it) =
O\((1+|t|)^{-m}\)$ for any $m\ge 0$.  Hence
$$\int_T^\infty\hat{g}(it)dt = O\((1+T)^{-m}\),~~T>0$$ and the
proposition follows because of our normalization
$\int_{{\R}}\hat{g}(it)dt=1.$\bx

\subsection{Partial Trace}

Let $C\ge 1$ be a fixed truncation parameter.  We will adjust it
only in the final step of our analysis.  Denote the usual
fundamental domain for $X$ as $$\F=\{z\in {\U} \mid |z+\bar{z}|\le
1 \le |z|\}$$ and its truncation as $$\F_{C}=\{z\in\F \mid
Im(z)\le C\}.$$
  Of course, $\F_C$ is compact and has area $\f{\pi}{3}-\f{1}{C}$;
  as $C\rightarrow\infty$, it exhausts $\F$.
If $X$ were compact, we would take the trace of the integral
operator $L_g$ by integrating the two expressions for $K(x,x)$ in
(\ref{specexpans}) over $\F$.  However, these integrals diverge
and we instead take a partial trace over $\F_C$:

$$\int_{{\F}}K(x,x)dx = \sum_{\g\in\bar{\G}} \int_{{\F}_C}g(x\i \g
x)dx$$
\begin{equation}
\label{firstparttrace} = \sum_{j\ge
0}\hat{g}(\nu_j)\int_{{\F}_C}|\phi_j(x)|^2dx +
\f{1}{4\pi}\int_{{\R}}\hat{g}(it)\int_{{\F}_C}|E_{1/2+it}(x)|^2dxdt.
\end{equation}

Up until now we have essentially followed Selberg's derivation of
the trace formula.  His next step is the understand the divergence
of each side of (\ref{firstparttrace}) as a function of $C$.  He
is able to cancel it from each side and derive a trace formula for
$\sum \hat{g}(\nu_j)$ by letting $C\rightarrow\infty$.  As a final
step, he uses a family of functions such as the $g_T$ to deduce
the spectral asymptotics.

This cancellation can be found  in great generality (e.g.
\cite{ArthurICM}), but the resulting formula is quite complicated.
We will now deviate from this path by keeping $C${\em fixed} for
now, and performing an analysis of (\ref{firstparttrace}) with the
$g_T$. By noting \begin{equation}\int_{{\F}_C}|\phi_j(x)|^2dx \le
\int_{{\F}}|\phi_j(x)|^2dx = 1,\label{sl2l2} \end{equation}
(\ref{firstparttrace}) becomes an inequality, giving a lower bound
on $\sum \hat{g}(\nu_j)$ (which, as we will make precise later,
roughly equals $N(T^2)$).

Since only a finite number of $\G$-translates of ${\F}_C$ neighbor
it, we may take $\sigma$ to be small enough so that $$g(x\i \g
x)\neq 0,~~x\in{\F}_C$$ only for the finitely-many $\g\in\bar{\G}$
which have a fixed-point on the boundary of ${\F}_C$.

\begin{lem}  For $\sigma$ small,

$$\left|   \sum_{\stackrel{\g\in\bar{\G}}{\g\neq I}}
\int_{{\F}_C}g(x\i \g x)dx \right| = O\( \sigma^2 \max{g}\).$$
\end{lem}
{\bf Proof:} By the above remark, it suffices to show
$$\int_{{\U}} |g(x\i \g x)|dx = O\( \sigma^2 \max{g} \)$$ for $\g$
a rotation.  Since rotations sweep distant points a proportional
amount, only $x$ of distance $O(\sigma)$ from the fixed point can
have $g(x\i \g x)\neq 0$.\bx

Applying (\ref{firstparttrace}), (\ref{sl2l2}), and
\propref{maxgt} we have that
\begin{eqnarray}
\int_{{\F}_C}g_T(x\i I x)\le & \sum_{j\ge 0}\hat{g}_T(\nu_j) &
\nonumber \\ & +
\f{1}{4\pi}\int_{{\R}}\hat{g}_T(it)\int_{{\F}_C}|E_{1/2+it}(x)|^2dxdt
& \nonumber \\ & + O(\sigma^2 g_T(I)). &
\end{eqnarray}

The next subsection is devoted to analyzing the Eisenstein series
term.  In \lemref{gl2eisbdlem} we show that it is $O(T\log T)$.
Thus for some $c_0>0$

$$\(\mbox{area}({\F}_C) - c_0\sigma^2\)g_T(I)
=\(\mbox{area}({\F}_C)
-
c_0\sigma^2\)\f{1}{4\pi}\int_{{\R}}\hat{g}_T(it)t\tanh(\pi t)dt $$
$$\le  \sum_{j\ge 0}\hat{g}_T(\nu_j) + O(T\log T). $$

{\bf Proof of \thmref{selbergsweylslaw}:}

 We need to show $$\lim_{T\rightarrow\infty}\f{N(T)}{T} =
\f{1}{12}.$$  The general theorem of \cite{Donnelly}, for example,
shows that \begin{equation}\limsup_{T\rightarrow\infty}
\f{N(T)}{T} \le \f{1}{12}.\label{sl2donn}\end{equation}

By \propref{gtandchar} \begin{eqnarray}
\int_{{\R}}\hat{g}_T(it)t\tanh(\pi t)dt = & \int_{-T}^T t\tanh(\pi
t)dt + O(T) & \nonumber \\ = & T^2+o(T^2). &
\end{eqnarray}
Again using (\ref{sl2donn}), $$\sum_{j\ge0}\hat{g}_T(\nu_j) \le
N((T+\sqrt{T})^2) + O(1).$$ This proves
$$\liminf_{T\rightarrow\infty}\f{N(T)}{T}\ge
\f{\f{\pi}{3}-\f{1}{C}-c_0\sigma^2}{4\pi}.$$ Taking
$\sigma\rightarrow 0$ and $C\rightarrow\infty$, we conclude the
Weyl law.\bx

\subsection{The Eisenstein Series Bounds}

Recall that the Eisenstein series $E_s(z)$ has the constant term
$$c(v,s)=\int_0^1 E_s(u+iv)du = v^s+\phi(s)v^{1-s},$$ where
$$Z(s)=\pi^{-s/2}\G(\f{s}{2})\zeta(s),$$ and
$$\phi(s)=\f{Z(2s-1)}{Z(2s)}~~=~~  \sqrt\pi\f{\G(s-\f
12)\zeta(2s-1)}{\G(s)\zeta(2s)}.$$

Define the {\em truncated Eisenstein series} for $u+iv\in\F$ as
\begin{equation}\L^CE_s(u+iv)=
\left\{\begin{array}{ll}  E_s(u+iv)-c(v,s), & v > C \\
                      E_s(u+iv), & v \le C.\end{array} \right.
\end{equation}
 The truncated $\L^CE_s(z)$ can be extended from $\F$ to $\G\backslash\U$;
  it decays rapidly as $v\rightarrow
\infty$, and so is in $L^2(\G\backslash\U)$.
  Thus for any $s_1,s_2\in\C$ which are not poles of the Eisenstein series $E_s(z)$,
$$\int_{\G\backslash\U}\L^CE_{s_1}(z)\L^CE_{s_2}(z)\frac{dudv}{v^2}
 < \infty.$$  In fact the same is true if we only
 truncate $E_{s_1}(z)$:
\begin{lem} With $s_1$ and $s_2$ as above,
 $$\int_{\G\backslash\U}\L^CE_{s_1}(z)\L^CE_{s_2}(z)\frac{dudv}{v^2}
=\int_{\G\backslash\U}\L^CE_{s_1}(z)E_{s_2}(z)\frac{dudv}{v^2}.$$
\end{lem}
{\bf Proof:} Unfolding, their difference
$$\int_{\G\backslash\U}\L^CE_{s_1}(z)\(E_{s_2}(z)
-\L^CE_{s_2}(z)\)\frac{dudv}{v^2}
=\int_C^\infty\frac{dv}{v^2}c(y,s_2)
\int_0^1\(E_{s_1}(z)-c(v,s_1)\)du$$
$$=\int_C^\infty\frac{dv}{v^2}c(v,s_2)\(c(v,s_1)-c(v,s_1)\)=0.$$
\bx

It will simplify notation to introduce $\d_C(z)$,
 the characteristic function of $\{z\in{\C}\mid \im z > C\}$.  Then for $z\in\F$
$$\L^CE_s(z)=E_s(z)-\d_C(z)c(v,s)
=\sum_{\g\in\G_{\infty}\backslash\G}\left[\im\!(\g z)^{s_1}
-\d_C(\g z)c(\im\!(\g z),s_1)\right].$$ Because the rightmost
expression is automorphic, it agrees with $\L^CE_s(z)$ for all
$z\in\U$; at most one of the terms indicated by the $\d_C$
functions
 is actually
subtracted.  We can thus unfold this series representation of
$\L^CE_{s_1}(z)$ in the integral:
$$\int_{\G\backslash\U}\L^CE_{s_1}(z)E_{s_2}(z)\frac{dudv}{v^2} =
 \int_{\G_{\infty}\backslash \U} E_{s_2}(z)\(v^{s_1}-
\d_C(z)c(v,s_1)\)\frac{dudv}{v^2}$$
$$=\int_0^{\infty}\frac{dv}{v^2}\(v^{s_1}-\d_C(z)c(v,s_1)\)
\int_0^1duE_{s_2}(z)$$
$$=\int_0^{\infty}\(v^{s_2}+\phi(s_2)v^{1-s_2}\)\(v^{s_1}
-\d_C(v)\(v^{s_1}+\phi(s_1)v^{1-s_1}\)\)\frac{dv}{v^2}$$
$$=\int_0^Cv^{s_2+s_1-2}dv + \phi(s_2)\int_0^C v^{s_1-s_2-1}dv $$
$$- \phi(s_1)\int_C^{\infty}v^{s_2-s_1+1}dv - \phi(s_1)\phi(s_2)
\int_C^{\infty}v^{-s_1-s_2}dv$$
\begin{equation}
=\frac{C^{s_2+s_1-1}}{s_2+s_1-1} +
\frac{\phi(s_2)C^{s_1-s_2}}{s_1-s_2} +
\frac{\phi(s_1)C^{s_2-s_1}}{s_2-s_1}
 + \phi(s_1)\phi(s_2)\frac{C^{1-s_2-s_1}}{1-s_2-s_1},
\label{origmaasssel}
\end{equation}
provided $\Re{s_1} > 1+\Re{s_2}$.  Each side has a meromorphic
continuation in $s_1,s_2$, and together are the {\em Maass-Selberg
relations}.

Because the Eisenstein series and hence their truncations are
holomorphic for $Re(s_1)=Re(s_2)=\f 12$, the Maass-Selberg
relations must have a removable singularity there. The limiting
value is
\begin{equation}
\int_{\G\backslash\U}|\L^CE_{1/2+it}(z)|^2\frac{dudv}{v^2}= 2\log
C-\f{\phi'}{\phi}(\f{1}{2}+it) + \f{C^{2it}\phi(\f
12-it)-C^{-2it}\phi(\f 12 + it)}{2it}. \label{gl2mslimitterm}
\end{equation}

The functional equations
$$Z(s)=Z(1-s)~,~Z(\bar{s})=\overline{Z(s)}$$ show that
\begin{equation}|\phi(\f 12 +it)|=
\left|\f{Z(it)}{Z(it+1)}\right|=1.\label{gl2unitary}\end{equation}
 Also,
\begin{equation}\f{\phi'}{\phi}(\f
12+it)=2\(\f{Z'}{Z}(2it)+\f{Z'}{Z}(-2it)\)\label{gl2scatlogdif}\end{equation}
and \begin{equation}\f{Z'}{Z}(s)=\sum_{\rho\mid
Z(\rho)=0}\f{1}{s-\rho} -\f 1s
-\f{1}{s-1},\label{melbrooks}\end{equation} a sum over the poles
and critical zeroes of $\zeta(s)$.

\begin{lem}
\label{gl2eisbdlem} $$\int_{{\R}}
\hat{g}_T(it)\int_{{\F}_C}|E_{1/2+it}(x)|^2 dx dt = O(T\log T).$$
\end{lem}
{\bf Proof:}  By the definition of the truncation,
$$\int_{{\F}_C}|E_{1/2+it}(x)|^2 dx \le
\int_{{\F}}|\L^CE_{1/2+it}(x)|^2 dx,$$  which is calculated in
(\ref{gl2mslimitterm}). In view of (\ref{gl2unitary}) and
(\ref{gl2scatlogdif}), it suffices to show $$\Re\int_{{\R}}
\hat{g}_T(it) \f{Z'}{Z}(2it)dt = O(T \log T).$$ This follows from
(\ref{melbrooks}) and the zero counting estimate (see
\propref{wind} or \cite{Titchmarsh})
$$\Re\int_{i(H-1)}^{i(H+1)}\sum_{\rho\mid Z(\rho)=0}\f{ds}{s-\rho}
= O(\#\{\rho \mid H-1\le \Im(\rho)\le H+1\}) = O(\log H).$$\bx


\section{Coordinates on $SL_3({\R})$ and $\H=SL_3({\R})/SO_3({\R})$}
\label{struct}

The Iwasawa decomposition $G=SL_3({\R})=NAK$ states that every
element $g\in G$ can be uniquely factored as $g=nak$, where $$n\in
N_0=\left\{\ttt1**01*001 \in G \right\},$$ $$a\in
A_0=\left\{\ttt{a_1}000{a_2}000{a_3}\in G~,~a_i>0\right\},$$ and
$$k\in K=SO_3({\R}).$$ The minimal standard parabolic is
$$P_0=N_0A_0=A_0N_0=\left\{\ttt***0**00* \in G \right\}.$$ There
are two associate maximal standard parabolics
$$P_1=\left\{\ttt******00* \in G \right\}~~,~~P_2=
\left\{\ttt***0**0** \in G \right\} .$$ Each can be decomposed
uniquely as $P_i=M_i'N_iA_i=N_iM_i'A_i$, where $$N_1=\left\{\ttt
10*01*001 \in G \right\},$$ $$M_1'=\left\{\ttt **0**000{\pm 1} \in
G \right\},$$ $$A_1=\left\{\ttt a000a000{a^{-2}} \in G
~,~a>0\right\},$$ $$N_2=\left\{\ttt 1**010001 \in G \right\},$$
$$M_2'=\left\{\ttt {\pm 1}000**0** \in G \right\},$$ and
$$A_2=\left\{\ttt {a^2}000{a\i}000{a\i} \in G ~,~a>0\right\}.$$

The {\em period} of an automorphic function $\psi\in
 L^2(\G\backslash G/K)$ along the parabolic $P$ is
defined as the integral \begin{equation} \psi_P(g)=\int_{\G\cap
N\backslash N} \psi(ng)dn.\label{perdef}\end{equation} If all
periods along each of the parabolics $P_0,P_1,$ and $P_2$ vanish,
then $\psi$ is called a {\em cusp} form.

%
%

  We shall use the following sets of coordinates for the Lie algebras.
  Let $${\bf a}_0=\left\{\ttt{h_1}000{h_2}000{h_3} \mid
h_1,h_2,h_3\in {\R}~,~ h_1+h_2+h_3=0\right\}$$ $$\simeq
\{H=(h_1,h_2,h_3)\in {\R}^3\mid h_1+h_2+h_3=0\}     .$$ It has
simple roots $\a_1,\a_2\in {\bf a}_0^*$, which are the following
linear functions on ${\bf a}_0$:
$$\a_1(H)=h_1-h_2~,~\a_2(H)=h_2-h_3.$$ There is also a third root
$\a_3=\a_1+\a_2$ which acts by $\a_3(H)=h_1-h_3$.  We can
similarly coordinatize

$${\bf a}_0^*=\left\{\l=(\ell_1,\ell_2,\ell_3)\in {\R}^3\mid
\ell_1+\ell_2+\ell_3=0      \right\}$$ and its complexification
$${\bf a}_{0{\C}}^*=\left\{\l=(\ell_1,\ell_2,\ell_3)\in {\C}^3\mid
\ell_1+\ell_2+\ell_3=0      \right\}.$$  Since these are subsets
of ${\C}^3$, we will use the standard norm there to define norms
in these spaces. The simple roots are bases of
 ${\bf a}_0^*$ and ${\bf a}_{0{\C}}^*$
as vector spaces over $\R$ and $\C$, respectively.  Also, they
have corresponding co-roots $$ \a_1^{\vee}=(1,-1,0)\in{\bf a}_0
~,~
 \a_2^{\vee}=(0,1,-1)\in{\bf a}_0.$$
%
%
%
The Weyl groups $\Omega({\bf a}_0)$ and $\Omega({\bf a}_0^*)$ are
isomorphic to the symmetric group $S_3$, and act in a compatible
way by permuting the standard basis vectors in ${\R}^3$ or
${\C}^3$ (see the chart in the appendix for more details).
 The
Cartan subgroups $A_1,A_2$ of the maximal parabolics have
one-dimensional Lie algebras, coordinatized by $${\bf
a_1}=\{H=(h,h,-2h)\mid h\in{\R}\} \simeq \{h\in\R\}$$ and $${\bf
a_2}=\{H=(2h,-h,-h)\mid h\in{\R}\} \simeq \{h\in\R\}.$$ These each
have one-dimensional dual spaces, ${\bf a_1}^*$ and ${\bf a_2}^*$
respectively.  Also, each dual has a special vector $\rho$, half
the sum of the positive roots:
$$\rho_0=\f{\a_1+\a_2+\a_3}{2}:(h_1,h_2,h_3)\mapsto h_1-h_3 ,$$
$$\rho_1 :(h,h,-2h)\mapsto 3h,$$ and $$\rho_2 :(2h,-h,-h)\mapsto
3h.$$ Under these coordinates $\rho_1$ and $\rho_2$ are naturally
identified with $\a_2$ and $\a_1$, respectively.

%
%
%

Each parabolic $P_i$ has the Langlands decomposition
$$P_i=N_iA_iM_i'$$ and logarithm maps  $H_i:G\rightarrow {\bf
a_i}$ such that \begin{equation}a=e^{H_i(a)}~~,~~a\in
A_i\label{logadef}\end{equation} and $$g\in N_ie^{H_i(g)}M_i'K.$$
These maps are well-defined despite the fact that the
decomposition $G=P_iK$ is in general not unique. Similarly for the
maximal parabolics, there are maps $m_1$ and $m_2$ mapping $G$
onto $M_1'/(K\cap M_1')$ and $M_2'/(K\cap M_2')$, respectively.

\subsection*{Truncation}

 Let $\D_P$ denote the set of simple roots which do not
vanish identically on $P$:
$$\D_{P_0}=\{\a_1,\a_2\}~,~\D_{P_1}=\{\a_2\}~,~
\D_{P_2}=\{\a_1\}.$$ The group $G$ can also be viewed as a
parabolic, with $\D_G=\{\}.$

We can now define Langlands' truncation (\cite{Langlands1966}, see
\cite{Arthur1980}).  For any parabolic $P= P_0,P_1,P_2$, or $G$,
define $\hat{\tau}_P(x)$ to be the characteristic function of
$$\{x=c_1\a_1^\vee+c_2\a_2^\vee\in
 {\bf a}_0 \mid c_i>0~\forall \a_i\in\D_P\}.$$
Thus, $\hat{\tau}_{P_0}$ is the characteristic function of
$$\{x=c_1\a_1^{\vee}+c_2\a_2^{\vee}\in {\bf a}_0 \mid
c_1,c_2>0\},$$ $$\hat{\tau}_{P_1}\mbox{ of }
\{x=c_1\a_1^{\vee}+c_2\a_2^{\vee}\in{\bf a}_0 \mid c_2>0\},$$
$$\hat{\tau}_{P_2}\mbox{ of }
\{x=c_1\a_1^{\vee}+c_2\a_2^{\vee}\in{\bf a}_0 \mid c_1>0\},$$ and
$$\hat{\tau}_{G}\mbox{ of } {\bf a}_0.$$

Let $C\in{\bf a}_0$ be a fixed parameter. The truncation of an
automorphic form $\psi$ is a sum over all standard parabolic
subgroups $$(\L^C\psi)(x):=\sum_P(-1)^{\dim A}\sum_{\g\in\G\cap
P\backslash\G}\hat{\tau}_P(H_0(\g x)-C)\int_{\G\cap N\backslash
N}\psi(n\g x)dn$$ \begin{equation}=\sum_P(-1)^{\dim
A}\sum_{\g\in\G\cap P\backslash\G}\hat{\tau}_P(H_0(\g
x)-C)\psi_P(\g x),\label{truncdef}\end{equation} which itself is
clearly an automorphic form. It can be proven that it also decays
rapidly in the cusp because of the way its constant terms have
been removed (see \cite{Arthur1980}).  Note that if $\psi$ is a
cusp form to begin with, by definition all its constant terms
$\psi_P$ vanish identically in proper parabolics, and thus
$\L^C\psi=\psi$.

%
%

\subsection*{Spectral Density}

We will later need to use the spectral density function $\b(\l)$
on $i{\bf a_0}^*$, which is a constant multiple of
 \begin{equation} \left| \f{
\G(\f 12+\l(\a_1^\vee)) }{\G(\l(\a_1^\vee)) } \f{ \G(\f
12+\l(\a_2^\vee)) }{\G(\l(\a_2^\vee)) } \f{ \G(\f
12+\l(\a_3^\vee)) }{\G(\l(\a_3^\vee)) }   \right|^2
.\label{spectraldensityform}\end{equation} This function will
appear below in (\ref{inver}); see \cite{Kolk} for more
information.  For $\l$ such that $\l(\a_1^\vee),\l(\a_2^\vee)$,
and $\l(\a_3^\vee)$ are all large, Stirling's formula shows that
$\b(\l)$ behaves as a constant times
$|\l(\a_1^\vee)\l(\a_2^\vee)\l(\a_3^\vee)|$.

The spectral density function is normalized so that
$$\int_{\|\l\|\le T}\b(\l)d\l\sim
 \f{T^{5/2}}{\G(7/2)(4\pi)^{5/2}}.$$  One can also view this as
normalizing the measure $d\l$.



%
%

\section{Convolution operators}
\label{convopsec}

 Every $g\in C_c^{\infty}(K\backslash G/K)$ acts
by convolution on $L^2(G/K)$: $$(L_gf)(x)=(f*g)(x)=\int_G f(y)
g(y\i x)dy=\int_{G/K} f(y) g(y\i x)dy.$$ The convolution operator
$L_g$ also acts on $f\in L^2(\G\backslash G/K)$ by
$$(L_gf)(x)=\int_{G/K} f(y) g(y\i x)dy$$ $$=\sum_{\g\in\G}
\int_{\G\backslash G/K} f(\g\i y) g(y\i \g x)dy=
\int_{\G\backslash G/K} f(y)K(x,y)dy,$$ where
$$K(x,y)=\sum_{\g\in\G}g(y\i \g x)$$ is the automorphic kernel.

Suppose furthermore that $g(x)$ is real; then $L_g$ is an operator
on $L^2(\G\backslash G/K)$ which commutes with the laplacian $\D$
and all other invariant differential operators.  We may thus
choose an orthonormal set $$\f{1}{\sqrt{vol(\G\backslash
G/K)}}=\phi_0, \phi_1,\phi_2,\ldots$$ of common eigenfunctions of
$L_g$ and the ring of invariant differential operators ${\cal D}$.
Continuing, we may resolve any multiplicities that remain by using
the Hecke operators, which commute with these other operators.

\begin{prop}(Selberg's Uniqueness Principle)(\cite{Selberg})
\label{selunique} If $\phi$ is a common eigenfunction of the ring
of invariant differential operators, then
$$(L_g\phi)(x)=\hat{g}(\phi)\phi(x),$$ where $\hat{g}(\phi)$ only
depends on $\phi$'s eigenvalues under ${\cal D}$.  In fact, one
can always find some $\l\in {\bf a}_{\C}^*$ such that the function
 $\phi_{\l}=e^{(\l+\rho)(H_0(x))}$
on $G/K$ has the same eigenvalues as $\phi$ under
any $D\in {\cal D}$.  This provides the formula
\begin{equation}
\hat{g}(\phi)=\hat{g}(\l)=(L_g\phi_\l)(I)=
\int_{G/K}g(x)e^{(\l+\rho)(H_0(x))}dx. \label{seluniq}
\end{equation}
\end{prop}
If $s$ is any permutation in the Weyl group $\Omega({\bf a}^*_0)$,
then in fact $\phi_\l(x)=e^{(\l+\rho)(H_0(x))}$ has the same
eigenvalues as $\phi_{s\l}$. The uniqueness principle thus implies
that $\hat{g}$ is invariant under $\Omega({\bf a}^*_0)$. Formula
(\ref{seluniq}) shows that the transform $\hat{g}(\l)$ is the
composition of an average over $N_0$, and a Fourier transform on
$A_0$.  More precisely, if $g$ is a function on $G/K$, denote

$$\bar{g}(a)=\int_{N_0} g(na)dn.$$  For a function $f$ on $A_0$,
define the Fourier transform as $$({\F}f)(\l) = \int_{A_0}
f(a)e^{\l(H_0(a))}da~~,~~\l\in i{\bf a^*_{0{\C}}}.$$ Then
$\hat{g}={\F}\(\bar{g}(\cdot) e^{-\rho(H_0(\cdot))}\)$.

We will make use of the following theorem in constructing our choices
of functions $g$.  Before stating it, let us introduce the notation
$$A(\sigma)=\{ a \in A_0 \mid \| H_0(a) \| \le \sigma \}.$$

\begin{thm}(\cite{Gangolli})
\label{gangolli} The map $g\mapsto \bar{g}$ provides a bijection
between the sets $$\{g \in C^{\infty}(K\backslash G/K) \mid supp~g
\subset KA(\sigma)K \}$$ and $$\{h \in C^{\infty}(A) \mid   supp~h
\subset A(\sigma)  \mbox{ and }  h(sa)=h(a)\mbox{ for all
}s\in\Omega({\bf a_0})\}.$$
\end{thm}

\subsection*{Construction of functions}

Let $\sigma>0$ be fixed for the remainder of this section.  It is
possible to find a non-negative function $g\in
C^{\infty}(K\backslash G/K)$ which is supported in $KA(\sigma)K$,
and whose transform $\hat{g}$ is non-negative on the joint
spectrum.  Having such a function without the positivity
requirement on $\hat{g}$ is straightforward; one then rescales and
convolves it with itself to achieve positivity (see lemma 6.2 of
\cite{Kolk} for example).  We shall again normalize $g$ so that
$$\int\!\!\!\!\int_{i{\bf a_0^*}} \hat{g}(\l)d\l=1.$$ Let us state
a property of $g$ which is implied by the positivity of $\hat{g}$
on $i{\bf a_0}^*$:

\begin{prop}
\label{mypos}  If $g\in C_C^{\infty}(K\backslash G/K)$ is such that
$\hat{g}(\l)\ge 0 $ for all $\l \in i{\bf a_0}^*$, then
$$\max_{x\in G} |g(x)| = g(I).$$
\end{prop}

\begin{remark}
{\em As the proof will demonstrate, the analogous fact is true for
the usual Euclidean Fourier transform.}
\end{remark}
{\bf Proof of \propref{mypos}:} The proof uses the inversion
formula
\begin{equation} g(x) =  \int_{i{\bf a_0}^*}
\hat{g}(\l) \overline{\tilde{\phi}_\l(x)}\b(\l)d\l, \label{inver}
\end{equation}

 where $$|\tilde{\phi}_\l(x)| \le
\tilde{\phi}_\l(I)=1$$ is a spherical function (e.g.
(\ref{spherfundef})) and $\b(\l)$ is the spectral density.
Trivially $$|g(x)| \le \int_{i{\bf a_0}^*} \hat{g}(\l)
\b(\l)d\l=g(I).$$ \bx

Now if $\Sigma$ is a measurable, bounded subset of
$i{\bf a_0}^*$ which is invariant under the Weyl group, define

\begin{equation} \hat{g}_{\Sigma}(\l)=\hat{g}*\chi_\Sigma(\l) = \int_\Sigma
\hat{g}(\l-\mu)d\mu. \label{defghatsig}
\end{equation}

 The function $\hat{g}_\Sigma(\l)$ is
roughly concentrated on $\Sigma$, especially for large $\Sigma$;
we shall use it to estimate the spectrum in $\Sigma$. Also
$\hat{g}_\Sigma(\l)$ decays rapidly, since $g_{\Sigma}$ is a
smooth function of compact support.
 If $\Sigma$ is open and its boundary has a finite Hausdorff length, then

\begin{equation}
\label{paley} \hat{g}_{t\Sigma}(\l) \ll (1+dist(\l,t\Sigma))^{-m}
\end{equation} (cf. p. 85 of \cite{Kolk}). Of course
$\hat{g}_{\Sigma}(\l)$ is the Fourier transform of
$$\bar{g}(a)\int_\Sigma e^{(-\mu-\rho)(H_0(a))}d\mu,$$ which is a
smooth function on $A_0$ whose support is contained in
$A(\sigma)$. It is furthermore invariant under the action of the
Weyl group $\Omega({\bf a_0})$ because $\Sigma$ is. Thus by
\thmref{gangolli} there exists a smooth, bi-$K$-invariant function
$g_\Sigma$ supported in $KA(\sigma)K$ such that
\begin{equation}
\label{defgsig}\overline{g_\Sigma}(a)=\bar{g}(a) \int_\Sigma
e^{-\mu(H_0(a))}d\mu
\end{equation} and
$$\widehat{g_\Sigma} = \hat{g}_\Sigma.$$ In summary, while not
changing the support of our function $g$, we can still smear its
transform $\hat{g}$ over $\Sigma.$  Of course $g_\Sigma$ becomes
more concentrated near the identity as $\Sigma$ gets larger;
 in the classical
Euclidean case this is analogous to multiplying a function by
Fejer's kernel.  In (\ref{bsig}) we will use a result analogous to
\propref{gtandchar} comparing $\hat{g}_\Sigma$ to $\chi_\Sigma$.

%
%

%
%

%
%


\section{The Partial Trace}
In addition to the discrete spectrum there is also a continuous
spectrum, furnished by Eisenstein series. Because it is
complicated we will just refer to its terms in the spectral
expansion as ``$Eis_g(x,y)$" until we need to be more
explicit:
\begin{equation}K(x,y)=\sum_{j\ge 0}
\hat{g}(\l_j)\phi_j(x)\phi_j(y) +Eis_g(x,y)=\sum_{\g\in\G}g(x\i\g
y).\label{gl3specexpn}\end{equation} Let $\F$ be a fundamental
domain for $\G\backslash G/K$ and $\F_C$ be a compact subset of
$\F$. Then
\begin{eqnarray}\int_{\F_C}K(x,x)dx = & \sum_{\g\in\G}\int_{\F_C} g(x\i\g
x)dx & \nonumber \\ = & \sum_{j\ge
0}\hat{g}(\l_j)\int_{\F_C}\phi_j(x)^2dx +\int_{\F_C}Eis_g(x,x)dx &
\end{eqnarray}
 Now,
$$\int_{\F_C}\phi_j(x)^2dx
\le
 \int_{\F}\phi_j(x)^2dx=1$$ so
\begin{equation}
\sum_{\g\in\G}\int_{\F_C}g(x\i \g x)dx \le \sum_{j\ge
0}\hat{g}(\l_j) + \int_{\F_C}Eis_g(x,x)dx. \label{firstineq}
\end{equation}

%
%

Now we will analyze the integrals $$\int_{\F_C}g(x\i \g x)dx $$
more systematically.  Note that we {\it do not} group these into
conjugacy classes as is done in the trace formula.

\begin{prop}
\label{mostvanish} For $\sigma$ sufficiently small and $supp~g
\subset A(\sigma)$, the integral
 $$\int_{\F_C}g(x\i \g x)dx =0$$ for
all but finitely many $\g$ -- those which have fixed points in the
closure of $\F_C$.
\end{prop}

\begin{prop}
\label{ellbd} If $\g\neq I$ has a fixed-point in the closure of
$\F_C$ and $supp~g \subset A(\sigma)$, then $$\int_{\F_C}g(x\i \g
x)dx \ll |\max g|vol(A(\sigma)).$$
\end{prop}

Although the implied constant above depends on $\F_C$, it may be
taken to be independent of $\g\in SL_3({\Z})$ since only finitely
many elements produce a non-zero integral.

Clearly $$\int_{\F_C}g(x\i I x)dx = vol(\F_C)\cdot g(I).$$ Now
take $g=g_\Sigma$ as in (\ref{defgsig}).  Using
Propositions~\ref{mypos} and \ref{ellbd}, the inequality
(\ref{firstineq}) becomes
\begin{equation}
\label{reintrosig}
 vol(\F_C) g_\Sigma(I) + O\(g_\Sigma(I)
vol(A(\sigma))\) \le \sum_{j=0}^\infty \hat{g}_\Sigma(\l_j) +
\int_{\F_C}Eis_{g_\Sigma}(x,x)dx.
\end{equation}
The proof of Theorem 8.5 in \cite{Kolk} shows that both

$$\left|     \sum_{j=0}^\infty \hat{g}_\Sigma(\l_j) - \# \{\im\l_j
\in \Sigma \}\right| \ll \int_{B\Sigma}\b(\l)d\l$$ and $$\left|
g_\Sigma(I) - \int_\Sigma \b(\l)d\l  \right| \ll
\int_{B\Sigma}\b(\l)d\l,$$ \begin{equation} \label{bsig} B\Sigma =
\{\l \in i{\bf a_0^*} \mid \mbox{dist}(\l,\partial \Sigma)\le 1
\},\end{equation} assuming $\Sigma$ is a bounded, open, and
Weyl-group invariant subset of $i{\bf a_0^*}$.  Now to control the
error term we shall assume $\Sigma$'s boundary is piece-wise
smooth (cf. Lemma 8.7 of \cite{Kolk}). Then
(\ref{spectraldensityform}) and (\ref{bsig}) show that for large
$t$ $$g_{t\Sigma}(I) \sim \int_{t\Sigma}\b(\l)d\l, $$ and by using
(\ref{reintrosig}), $$\(\int_{t\Sigma}\b(\l)d\l\)\(vol(\F_C)+
O(vol(A(\sigma)))\)(1+o(1))~~~~~~~~~~~~~~~~~~~~~~~~~$$
$$~~~~~~~~~~~~~~~\le \#\{\im\l_j\in t\Sigma\}+ \int_{\F_C}
Eis_{g_{t\Sigma}}(x,x)dx.$$ We will see at the end of
\secref{eisbdsec} that
\begin{equation}
\label{eisbound} \int_{\F_C} Eis_{g_{t\Sigma}}(x,x)dx =
o\(\int_{t\Sigma} \b(\l)d\l \).
\end{equation}
{\bf Proof of \thmref{mainresult}:} Taking $t\rightarrow\infty$
$$\liminf_{t\rightarrow\infty} \f{\#\{\im\l_j \in
t\Sigma\}}{\int_{t\Sigma}\b(\l)d\l} \ge
vol(\F_C)+O(vol(A(\sigma))).$$ We had insisted that $supp~g\subset
A(\sigma)$, so taking $\sigma\rightarrow 0$ and exhausting $\F$
through compact sets $\F_C$ we conclude that
\begin{equation}
\label{speclb}vol(\F)\le \liminf_{t\rightarrow\infty}
\f{\#\{\im\l_j \in
t\Sigma\}}{\int_{t\Sigma}\b(\l)d\l}.\end{equation}
 If $$\Sigma =
\{\l=(\ell_1,\ell_2,\ell_3) \in i{\bf a_0^*} \mid \|\l\|^2
 \le 1\},$$ then this indicates
$$vol(\F) \le \liminf_{T\rightarrow\infty}
\f{N(T)}{\(\f{T}{4\pi}\)^{5/2}\f{1}{\G(7/2)}}$$ because the
Laplace eigenvalue of $e^{(\l+\rho)(H(g))}$ is
$1-\f{\ell_1^2+\ell_2^2+\ell_3^2}{2}$.\footnote{A check of this
normalization is provided by the compact case, where there are no
Eisenstein series and the Weyl law is already known.} (Note that
$|\re{\ell_j}|< \f{1}{2}$ by unitary -- \cite{JacSha}.)

 The upper bound \begin{equation}\limsup_{T\rightarrow\infty}
\f{N(T)}{\(\f{T}{4\pi}\)^{5/2}\f{1}{\G(7/2)}} \le
vol(\F)\label{gl3donnupbd}\end{equation} due to \cite{Donnelly}
shows in fact that \begin{equation}N(T)\sim
\(\f{T}{4\pi}\)^{5/2}\f{vol(\F)}{\G(7/2)}.\label{gl3pfweyllawstat}\end{equation}
\bx




\begin{thm} (Spectral Equidistribution)
\label{genct} If $\Sigma$ is a bounded, open, Weyl-group invariant
subset of $i{\bf a_0^*}$ with a piece-wise smooth boundary, then
\begin{equation}\lim_{t\rightarrow\infty} \f{\#\{\im\l_j \in
t\Sigma\}}{\int_{t\Sigma}\b(\l)d\l}=
vol(\F).\label{equistate}\end{equation}
\end{thm}
{\bf Proof:} The lower bound is in (\ref{speclb}).  To prove the
upper bound, we will use this lower bound along with the
asymptotics of the Weyl law for counting in balls.   Let $r$ be
the length of the longest vector in $\Sigma$ (which is bounded by
assumption). By defining $\Sigma_c$ as the complement of $\Sigma$
inside the $r$-ball $B_r$, we have that $$\#\{\im\l_j \in
t\Sigma\} + \#\{\im\l_j \in t\Sigma_c\} = \#\{\im\l_j \in
t(\Sigma\cup \Sigma_c) \}.$$  Given any $\e>0$ we can find $t$
large enough so that \begin{equation}\#\{\im\l_j \in t\Sigma\} \le
(1+\e)\(\int_{t(\Sigma\cup \Sigma_c)}\b(\l)d\l -
\int_{t\Sigma_c}\b(\l)d\l\),\label{masscons}\end{equation} which
implies the upper bound when $\e\rightarrow 0$.\bx

 {\bf Proofs of Theorems \ref{nonselfdual} and
\ref{aeram}}: If the cusp form $\phi_j$ is not tempered, then its
spectral parameter $\l=(\ell_1,\ell_2,\ell_3)$ is not a
purely-imaginary vector.
 By the classification of the unitary dual, we have equality
of the sets $$\{\ell_1,\ell_2,\ell_3\} = \{-\bar{\ell}_1,
-\bar{\ell}_2,-\bar{\ell}_3\}.$$  We know that $|Re(\ell_i)|<\f
12$ by unitary (\cite{JacSha}), so the vectors $\l$ all lie near
the hyperplanes defined by
$$\l(\a_1^\vee)=0~,~\l(\a_2^\vee)=0~,~\l(\a_1^\vee+\a_2^\vee)=0.$$

Similarly, if $\phi$ is self dual then $$\{\ell_1,\ell_2,\ell_3\}
= \{\mu,0,-\mu\} \mbox{ for some }\mu.$$ This again constrains
$\l$ to lie along a hyperplane.

However, by taking the shape $\Sigma$ in \thmref{genct} to be
further and further away from any fixed hyperplane, we conclude
that no hyperplane has a positive percentage of the spectrum near
it. Thus each of the exceptional sets we are considering is of
measure zero compared to the rest of the spectrum.\bx

A related argument was used in \cite{Kolk} for cocompact
 subgroups $\G$.

{\bf Proof of Theorem \ref{native}}: All Gelbart-Jacquet lifts are
self-dual forms on $SL_3$.

 Another proof would be simply by
counting: the lift quadruples the Laplace eigenvalue of a cusp
form on $SL_2({\Z})\backslash \U$, and the number of these with
$SL_3({\R})$-Laplace eigenvalue $\le T$ is $O(T)$ by Selberg's
theorem \ref{selweyl}.  Though forms $\psi$
 on congruence
covers of $SL_2({\Z})\backslash \U$ may also lift to
$SL_3({\Z})\backslash SL_3({\R})/SO_3({\R})$, such $\psi$ are
actually twists of forms on $SL_2({\Z})\backslash \U$, so one need
only count the lifts from $SL_2({\Z})\backslash \U$ itself, and
not from other congruence covers.
  \bx

\section{Eisenstein Series}

These periodized functions on $\G\backslash G/K$ are constructed
using the Langlands decompositions of $G$'s parabolics. Recall
that any element $p\in P=NAM'$ factors uniquely into
 a product $p=nam$ of elements from their respective subgroups.
Writing the diagonal matrices $a=e^{H_P(p)}$, we get a map which
extends to $G$: $$g \mapsto H_P(g) ~,~g\in N e^{H_P(g)} M' K.$$
For the parabolics $P_1$ and $P_2$, where $M'\simeq GL_2({\R})$,
there are corresponding maps $m:G\rightarrow M'/(K\cap M')$ (see
(\ref{logadef})).

If $\l\in{\bf a}_{0{\C}}^*$ and $g\in G$, the minimal parabolic
Eisenstein series is defined as
\begin{equation}
E(P_0,g,\l)=E(P_0,g,\l,1)=\sum_{\G\cap P_0\backslash \G}
e^{(\l+\rho_0)(H_0(\g g))}.\label{minparaeisdef} \end{equation}
This sum only converges when $\a_1(\l)$ and $\a_2(\l)$ have large
real parts, but it has a meromorphic continuation to all of ${\bf
a}_{0{\C}}^*$.

Since $$\G\cap M_1'\backslash M_1'/(K \cap M_1')\simeq
GL_2({\Z})\backslash \U,$$ the discrete eigenfunctions on the
former are just the even\footnote{i.e. $f(x+iy)=f(-x+iy)$.}
discrete eigenfunctions for $SL_2({\Z})$.  Take such a cusp form
and $\l\in{\bf a}_{1{\C}}^*$, and define the maximal parabolic
Eisenstein series \begin{equation}E(P_1,g,\l,\phi)=\sum_{\G\cap
P_1\backslash \G}
 e^{(\l+\rho_1)(H_1(\g g))}\phi(m_1(\g g)).\label{maxparaeisdef}\end{equation}
  There is a similar
Eisenstein series $E(P_2,g,\l_1,\phi)$ for the other maximal
parabolic $P_2$, related through a functional equation.  Each
series again is initially only defined for certain values of $\l$
but extends via a meromorphic continuation to ${\bf a}_{\C}^*$
(\cite{Langlands1966},\cite{Langlands1976}).

We are now in a position to define what ``$Eis_g(x,y)$" is.  The
spectral expansion of the automorphic kernel is $$ K(x,y)=
\sum_{\stackrel{\phi_j\mbox{ an $L^2$ discrete
eigenfunction}}{\mbox{on } SL_3({\Z})\backslash
SL_3({\R})/SO_3({\R})}}
\hat{g}(\l_j)\phi_j(x)\phi_j(y)+Eis_g(x,y), $$ and up to a
normalizing constant for the measure, $Eis_g(x,y)$ is
(\cite{ArthurICM})

$$ \f{1}{3(2\pi i)^2} \int\!\!\!\!\int_{i{\bf a_0^*}} \hat{g}(\l)
E(P_0,x,\l,1)\overline{E(P_0,y,\l,1)}d\l   +
~~~~~~~~~~~~~~~~~~~~~~~~~~~~~~~~~~~~~~~~~~~~~$$ \begin{equation}
\f{1}{2\pi i}\sum_{\stackrel{\phi_j\mbox{ an even $L^2$ discrete
eigenfunction}}{\mbox{on } SL_2({\Z})\backslash {\U}
~,~\D\phi_j=(\f14-\nu_j^2)\phi_j}} \int_{i{\bf
a}_1^*}\hat{g}(\l+(\nu_j,-\nu_j,0)) E(P_1,x,\l,\phi_j)
\overline{E(P_1,y,\l,\phi_j) }d\l. \label{eisdef}
\end{equation}
We also may assume that each $\phi_j$ on $\U$ is a Hecke
eigenform.  There is a beautiful formula for the inner-products of
truncated Eisenstein series, due to Langlands. It generalizes the
Maass-Selberg formula (\ref{origmaasssel}) for $SL_2({\R})$. See
\cite{Arthur1980} for details.

\begin{thm}Langlands' inner product formula (\cite{Langlands1966},
\cite{Arthur1980})
$$\int_{\G\backslash
G/K}\left[(\L^CE)(P,x,\phi,\l_1) \right]\left[(\L^CE)(P,x,\phi,\l_2)
\right]dx$$ $$=\sum_{P\sim P'\mbox{ associate}}
vol({\bf a}'/<\a^{\vee}\mid \a \in\D_{P'}>)\times$$
 $$\sum_{s_1,s_2\in\Omega({\bf
a},{\bf a}')}\f{e^{(s_1\l_1+s_2\l_2)(C)}}{\prod_{\a\in
\D_{P'}}(s_1\l_1+s_2\l_2)(\a^{\vee})}\left<M(s_1,\l_1)\phi,
M(s_2,\l_2)\phi\right>.$$
\end{thm}
 Here $M(s,\l)$ is an intertwining operator,
which sends $\phi$ to a cusp form on the potentially-different
parabolic $P'$.  We will discuss it in the instances it arises for
us, where it essentially acts as scalar multiplication.  The last
expression $\left<\psi,\psi'\right>$
 is an inner product over $\G\cap
M'\backslash M'$, and the Weyl group $\Omega({\bf a},{\bf a}')$ is
the set of isomorphisms of ${\bf a}$ to ${\bf a}'$ coming from
restrictions of elements in $\Omega({\bf a}_0)$.

\subsection*{Interlude: $SL_2$}
Write $C=(c,-c), c>0$, and
let $\l_1=(it+\e,-it-\e),\l_2=(-it,it)$ so that
$$\int_{\F}\L^C E(P,g,\l_1,\phi)E(P,g,\l_2,\phi)dg$$
is a constant times
$$\f{e^{(\e,-\e)C}}{2\e} +
 \f{e^{(2it+\e,-2it-\e)C}}{4it+2\e}R(-2it) +
 \f{e^{(-2it-\e,2it+\e)C}}{-4it-2\e}R(2it+2\e) +
\f{e^{(-\e,\e)C}}{-2\e}R(-2it)R(2it+2\e).$$

Here the intertwining operator is
\begin{equation}
R(s)=\sqrt{\pi}\f{\G(\f{s}{2})}{\G(\f{s+1}{2})}
\f{\zeta(s)}{\zeta(s+1)}=\f{Z(s)}{Z(s+1)}
\label{ratiodef}
\end{equation}
with
$$Z(s)=\pi^{-s/2}\G(\f s2)\zeta(s)=Z(1-s).$$

Take $\e\rightarrow 0$ so that $\l_1\rightarrow \bar{\l}_2$. Then
the last expression approaches
$$\int_{\F}|\L^CE(P,g,\l)|^2dg=(const) \left[4c-\f{R'}{R}(2it) +
\f{e^{4itc}}{4it}R(-2it) - \f{e^{-4itc}}{4it}R(2it)\right].$$ In
(\ref{origmaasssel})
 we derived this from a direct calculation.
It is a key step in Selberg's trace formula for $SL_2({\Z})$.

\section{Bounding the Eisenstein Contribution}
\label{eisbdsec}

In this final section we will complete the proofs by establishing
the estimate in (\ref{eisbound}).  From (\ref{eisdef}) it is
sufficient to show that both
\begin{equation}
\int_{\F_C} \int\!\!\!\!\int_{i{\bf a_0^*}}
\hat{g}_{t\Sigma}(\l)|E(P_0,x,\l,1)|^2d\l dx  =
o\(\int_{t\Sigma}\b(\l)d\l  \) ,\label{minpbd}
\end{equation} and
\begin{equation}
 \int_{\F_C} \int_{i{\bf a_1^*}} \sum_{\stackrel{\stackrel{\phi_j\mbox{
 an even $L^2$
discrete}}{\mbox{ eigenfunction}}}{\stackrel{\mbox{on }
SL_2({\Z})\backslash {\U}}{\D\phi_j=(\f14-\nu_j^2)\phi_j}}}
\hat{g}_{t\Sigma}(\l+(\nu_j,-\nu_j,0))|E(P_1,x,\l,\phi_j)|^2d\l dx
=o\(\int_{t\Sigma}\b(\l)d\l  \) \label{maxpbd}.
\end{equation}

To do this we shall use the rapid decay of $\hat{g}_{t\Sigma}$ in
(\ref{paley}).  Among other benefits, this allows us to
interchange the order of integration.  Keep in mind that
$\hat{g}_{t\Sigma}$ is roughly the characteristic function of
$t\Sigma$.

\subsection*{Some background on L-functions}

We will require some information about the density of zeroes of
certain L-functions: the Riemann $\zeta$ function and the standard
L-function of an even cusp form $\phi$ on $SL_2({\Z})\backslash
\U$, $L(s,\phi)$.  Each is defined as an Euler
product\footnote{Recall that $\phi$ is assumed to be both a Hecke
and Laplace eigenfunction.} over the primes
\begin{equation}
\zeta(s)=\prod_p (1-p^{-s})\i~~,~~L(s,\phi)=\prod_p (1-\a_p
p^{-s})\i(1-\a_p\i p^{-s})\i, \label{ldef} \end{equation} where
the $\a_p$ satisfy the bound $|\a_p| \le p^{5/28}$ (\cite{BDHI}).
Each L-function can be completed:
$$Z(s)=\G_{{\R}}(s)\zeta(s)~,~\L(s,\phi)=
\G_{{\R}}(s+\nu)\G_{{\R}}(s-\nu)L(s,\phi),$$ where
$$\G_{{\R}}(s)=\pi^{-s/2}\G(s/2)$$ and $\nu$ is related to
$\phi$'s Laplace eigenvalue by $$\D\phi=(\f14-\nu^2)\phi.$$ With
this convention $Z(s)=Z(1-s), \L(s,\phi)=\L(1-s,\phi)$, and each
is entire except for the simple poles of $Z(s)$ at $s=0,1$.

The following estimate will be used to bound Eisenstein series
integrals later.  The analogous statement for $Z(s)$ is classical
(e.g. see \cite{Titchmarsh}) and the proof for $\L(s,\phi)$ is
essentially identical. However we do not know of a reference in
the literature and include it for completeness.

\begin{prop}\label{wind}
For $T\in\R$
\begin{equation}\Re\int_{T-1}^{T+1}\f{\L'}{\L}(1+it,\phi)dt \ll
\log(|T|+|\nu|).\label{propwindstat}\end{equation}
\end{prop}
{\bf Proof}: Using entirety of $\L(s,\phi)$ and its Mittag-Leffler
expansion, $$\f{\L'}{\L}(s)=\f{\G_{{\R}}'}{\G_{{\R}}}(s+\nu)
+\f{\G_{{\R}}'}{\G_{{\R}}}(s-\nu) +\f{L'}{L}(s,\phi) =
\sum_{\{\rho \mid \L(\rho)=0 \}}\f{1}{s-\rho}$$ (the sum of the
zeroes is actually only conditionally convergent, so the term with
$\rho$ should always be summed with the term containing the zero
at $1-\rho$).  From the Euler product, $$\f{L'}{L}(s,\phi)= -
\sum_p\sum_{n=1}^{\infty}(\a_p^n + \a_p^{-n})p^{-ns}\log p,$$ so
the bound $|\a_p|\le p^{5/28}$ implies that $|\f{L'}{L}(2+i
t,\phi)|$ is uniformly bounded in both $\nu$ and $t$.
 By Stirling's formula,
$$\f{\G_{{\R}}'}{\G_{{\R}}}(s)=-\frac{1}{2}\log \pi +
\frac{1}{2}\log s + O(1/|s|)$$ and
$$\f{\G_{{\R}}'}{\G_{{\R}}}(2+it+\nu)+\f{\G_{{\R}}'}{\G_{{\R}}}(2+it-\nu)
\le \log(|t|+|\nu|) + O(1).$$ Thus,
\begin{equation}
\sum_{\rho}\f{1}{1+|t-\g|^2} \ll \sum_{\rho}\f{1}{2+it-\rho} \le
\log(|t|+|\nu|) +O(1). \label{zerosumbd}
\end{equation}
 It follows that there are no more than
$O(\log(|T|+|\nu|))$ zeroes between with imaginary part between
$T-1$ and $T+1$. By writing
$$\Re\int_{T-1}^{T+1}\f{\L'}{\L}(1+it,\phi)dt =
\Re\int_{T-1}^{T+1} \sum_{\rho}\f{1}{1+it-\rho}dt $$
$$=\Re\int_{T-1}^{T+1} \sum_{|\rho-iT|\le 2}\f{1}{1+it-\rho}dt
+\Re\int_{T-1}^{T+1} \sum_{|\rho-iT| > 2}\f{1}{1+it-\rho}dt,$$
invoking (\ref{zerosumbd}), and using the fact that
$$\re\int_{1+i(T-1)}^{1+i(T+1)}\f{ds}{s-\rho} =O(1),$$ we bound
each of the terms above by $O(\log(|T|+|\nu|))$.\bx

\subsection*{Minimal Parabolic Eisenstein Series}

We have appended a table of the 36 terms in the Langlands
inner-product formula (also referred to as the Maass-Selberg
relations) for the minimal parabolic.

The calculation is aided by the identity
\begin{equation}
M(s,(\ell_1,\ell_2,\ell_3)) =\prod_{\stackrel{1\le i<j\le
3}{s(i)>s(j)}} R(\ell_i-\ell_j). \label{scatprod}
\end{equation}  This is a special case of a general result of Langlands
(\cite{Euler},pp. 36-47;\cite{Bluelands}, p.
134;\cite{ArthurICM},p. 854);
 the function $R$ here
is the same as the one used above in (\ref{ratiodef}).

Using a limiting procedure as in the interlude we can compute:
\begin{prop}(Diagonal terms)
\label{diagterms} Let $C=(c,0,-c),
\l_1=(it_1+\e_1,it_2+\e_2,it_3+\e_3), \l_2=(-it_1,-it_2,-it_3)$.
Then $$\lim_{\e_1,\e_2,\e_3\rightarrow 0} \sum_{s\in\Omega({\bf
a}_0)}
\f{e^{s(\l_1+\bar{\l}_2)(C)}}{[s(\l_1+\bar{\l}_2)(\a_1^{\vee})]
[s(\l_1+\bar{\l}_2)(\a_2^{\vee})]} \left< M(s,\l_1),
M(s,\l_2)\right>$$ $$= 3 c^2-2 c \f{R'}{R}(it_1-it_2) -2 c
\f{R'}{R}(it_2-it_3) -2 c \f{R'}{R}(it_1-it_3)$$
$$+\f{R'}{R}(it_1-it_2)\f{R'}{R}(it_2-it_3)$$

$$+\f{R'}{R}(it_1-it_3)
\f{R'}{R}(it_2-it_3)
+\f{R'}{R}(it_1-it_2)\f{R'}{R}(it_1-it_3).$$
\end{prop}

\begin{prop}
\label{ptwiseminbd} For any $\e>0$ and $\l\in i{\bf a_0}^*$ with
$\|\l\|$ large, we have that $$\int_{\F} |\L^CE(P_0,x,\l)|^2dx =
O_{\e}(\|\l\|^\e).$$
\end{prop}
{\bf Proof:} We start with some estimates on $\zeta(s)$.  The
following may be found in Chapter 3 of \cite{Titchmarsh}: There is
an absolute constant $\kappa>0$ such that
\begin{equation}
\f{1}{\log(|t|+2)} \ll |\zeta(\sigma+it)| \ll \log(|t|+2)
\label{zetaupdown} \end{equation}
 and
$$\left|\f{\zeta'}{\zeta}(\sigma+it)\right| \ll \log(|t|+2)$$ in
the region $\sigma \ge 1-\f{\kappa}{\log(|t|+2)}$. Thus
$$R(s)=\f{Z(1-s)}{Z(1+s)} =
\f{\pi^{-(1-s)/2}\G\(\f{1-s}{2}\)\zeta(1-
s)}{\pi^{-(1+s)/2}\G\(\f{1+s}{2}\)\zeta(1+s)}$$ $$\ll
\(1+|\im\!{s}|\)^{-\re\!{s}}\log(|\im\!{s}|+2)^2$$ and
$$\f{R'}{R}(s) \ll \log(|\im\!{s}|+2)$$ for $$|\re\!{s}| \le
\f{\kappa}{\log(|\im\!{s}|+2)}.$$  Of course both $R(s)$ and
$\f{R'}{R}(s)$ are analytic in this region because of the
nonvanishing.

Fix $\l=(it_1,it_2,it_3)$.  Then each of the six terms from
\propref{diagterms} is trivially $O_{\e}(\|\l\|^\e)$.  Of the
remaining thirty terms (see the chart in the appendix), some may
have singularities when either of the denominators
$$s(\l_1+\bar{\l}_2)(\a_1^\vee)~~\mbox{or}~~s(\l_1+\bar{\l}_2)(\a_2^\vee)$$
vanish.  Nevertheless, the sum of all of the terms represents a
holomorphic function for all $\l \in i{\bf a_0}^*$, so the poles
cancel with other terms. If $\l$ is such that each
$$|s(\l_1+\bar{\l}_2)(\a_j^\vee)| \ge \f{\kappa}{2
\log(\|\l\|+2)},$$ then each term is trivially $O_{\e}(\|\l\|^\e)$
as well (the numerators have modulus one when $\re~{\l}=0$).

Otherwise, if some denominator is small, we will take estimates
further away and appeal to the maximum principle.  Take a small
neighborhood in $$s_1=s(\l_1+\bar{\l}_2)(\a_1^\vee)~~,~~
s_2=s(\l_1+\bar{\l}_2)(\a_2^\vee).$$   Now, suppose that $$|s_1|
\le \f{\kappa}{2\log(\|\l\|+2)} \le |s_2|$$ (the other cases have
almost-identical proofs). Then for
$|s_1|=\f{\kappa}{2\log(\|\l\|+2)}$, (\ref{zetaupdown}) shows that
each term is $O_{\e}(\|\l\|^\e)$.  The maximum modulus principle
shows that this bound holds uniformly for $|s_1| \ll \log \|\l\|$.
\bx

\subsection*{Maximal Parabolic $P_1$}

Now we turn to the sum in (\ref{maxpbd}).  We will first consider
the simpler case that $\phi$ is a cusp form.  If $s$ is the lone
permutation in $\Omega({\bf a}_1,{\bf a}_2)$, then the
intertwining operator acts as $$M(s,\l)\phi =
R(\l(\a^\vee),\phi)\phi',$$ where $\phi$ and $\phi'$ are cusp
forms on $M_1$ and $M_2$ coming from the same even cusp form on
$SL_2({\Z})\backslash \U$, and
$$R(s,\phi)=\f{\L(s,\phi)}{\L(s+1,\phi)}$$ is a ratio formed from
$\phi$'s completed standard L-function (\ref{ldef}).  As before
with (\ref{scatprod}), since $\phi$ is everywhere unramified, this
can be derived from Langlands' formula
(\cite{Euler},\cite{Bluelands}).  Then the Maass-Selberg relations
take the following form:

\begin{prop} The inner product
$$\int_{\F}\left| \L^CE(P,g,\l,\phi) \right|^2dg$$
is a constant multiple of $3/2 c - \f{R'}{R}(\l(\a^\vee),\phi)$.
\end{prop}
{\bf Proof}:
The inner product formula is a constant times
$$\lim_{\e\rightarrow 0} \f{e^{(2it+2\e)c}}{2it+2\e} -
\f{e^{-(2it+2\e)c}}{2it+2\e}\f{R(2it+2\e)}{R(2it)} .$$
\bx

\begin{prop} The integral

$$\int_{t=T-1}^{T+1}\int_{\F}\left| \L^CE(P,g,\l,\phi_j)
\right|^2dgdt = O(log(T+\|\l_j\|)).$$ \label{maxcuspbd}
\end{prop}
{\bf Proof}: Using the previous proposition, this requires only
the estimate on $\int_{T-1}^{T+1} \f{\L'}{\L}(1+it,\phi)dt$ in
\propref{wind} (cf. (\ref{gl2scatlogdif})).
 \bx


\begin{remark}\label{specsum} {\em (On summing over the eigenvalues) The formula
(\ref{eisdef}) includes a sum over the $SL_2({\Z})\backslash \U$
spectrum as well. Through a more-precise statement of Selberg's
\thmref{selweyl} such as $$N(T)=\f{1}{12}T^2+O(T\log T),$$ we can
bound the spectral points $\nu_j$ in the interval $[\nu-1,\nu+1]$
by $O(\sqrt{|\nu|}\log|\nu|)$.}
\end{remark}

\subsection*{The Constant Function on $SL_2$}

There is only one Eisenstein series left to estimate,
$$E(P_1,g,\l,1),$$ the Eisenstein series induced from the constant
function on the maximal parabolic $P_1$.

The constant functions on $SL_2({\Z})\backslash \U$ may be viewed
as multiples of $$Res_{s=1} E_s(z) =\f 12
Res_{s=1}\sum_{(c,d)=1}\f{y^s}{|cz+d|^{2s}}.$$ We may ignore the
actual value of the constant since we are only trying to get an
order-of-magnitude estimate. Thus, $E(P_1,g,\l,1)$is a constant
multiple of the residue \begin{equation}Res_{\d=0}E(P_0,g,(\f
12+it+i\d,-\f 12 + it -i\d,-2it),1).\label{constres}\end{equation}

Taking the residue of the inner product
$$\int_{\F}|\L^CE(P_1,g,\l,1)|^2dg$$ is more complicated, though a
limiting value must exist since the Eisenstein series is
meromorphic there. We will explicitly see this cancellation
occurring. Slicker arguments are possible, but we will present a
detailed proof for the sake of clarity.

Let $$\l_1=(\f 12+it+i\d_1+\e,-\f 12+it-i\d_1+\e,-2it-2\e),$$ and
$$\l_2=(\f 12-it+i\d_2~,~-\f 12-it-i\d_2~,~2it).$$ Then we are
interested in $$\lim_{\d_1,\d_2\rightarrow 0}
\d_1\d_2\sum_{s_1,s_2\in\Omega({\bf a}_0)}
 \f{e^{(s_1\l_1+s_2\l_2)(C)}}{(s_1\l_1+s_2\l_2)(\a_1^{\vee})
(s_1\l_1+s_2\l_2)(\a_2^{\vee}) }\left<M(s_1,\l_1)
M(s_2,\l_2)\right>.$$ This expression is holomorphic in $\d_1$ and
$\d_2$ near zero, so it does not matter how we take the limits
$\d_1,\d_2\rightarrow 0$. We will first take the residue in
$\d_1$. Of the 36 terms, we will of course ignore those which do
not have a pole at $\d_1=0$. These can occur only in
$M(s_1,\l_1)$, for the denominators do not yet vanish (note
$\e,\d_2\neq0$ at this stage).  Since $$M(s_1,\l_1)=
\prod_{\stackrel{i<j}{s_1(i)>s_1(j)}}R(\l_{1_i}-\l_{1_j}),$$ these
permutations $s_1$ must interchange 1 and 2, which forces
$$s_1\in\{(12),(13),(321)\}.$$

Next, when the residue at $\d_2=0$ is taken, poles can occur in
two different ways: $M(s_2,\l_2)$ might have a pole, or a
denominator might vanish. The former occurs for
$s_2\in\{(12),(13),(321)\}$ as well.

\begin{prop} Let $s_1\in\{(12),(13),(321)\}$.  If one of the
denominator terms $$(s_1\l_1+s_2\l_2)(\a_1^{\vee})=0$$ or
$$(s_1\l_1+s_2\l_2)(\a_2^{\vee})=0,$$ then
$$s_2\in\{e,(23),(123)\}=S_3 - \{(12),(13),(321)\},$$ i.e.
$M(s_2,\l_2)$ has no pole at $\d_2=0$.
\end{prop}
{\bf Proof}: Now that
$$\l_1=(\f 12+it +\e,-\f 12+it+\e,-2it-2\e)$$
and
$$\l_2=(\f 12-it,-\f 12-it,2it),$$
the only possible way to get consecutive entries of
$s_1\l_1+s_2\l_2$ to equal is if $$s_1=s_2(12).$$
\bx

We see that only the following terms have poles:

$$\mbox{two from }M~'s: s_1,s_2\in\{(12),(13),(321)\}$$ or
$$\mbox{one from }M,\mbox{ one from a denominator}:
s_1=(12),s_2=e~~~~,~~~~ s_1=(13),s_2=(123).$$ Note that
$(321)\times(23)$ fails to have poles in the denominators. The 11
terms (see the chart in the appendix) are:

$$(12)\times(12):
\f{e^{(-1+\e,1+\e,-2\e)(c,0,-c)}}{(-2)(1+3\e)}$$

$$(12)\times(13):
\f{e^{(-\f 12+3it+\e,\e,\f 12-3it-2\e)(c,0,-c)}}{(-\f 12+3it)(-\f 12
+3it+3\e)} R(\f 12-3it) R(-\f 12-3it)$$

$$(12)\times(321):
\f{e^{(-1+\e,\f 12+3it+\e,\f 12-3it-2\e)(c,0,-c)}}{
(-\f 32 -3it)(6it+3\e)}R(\f 12-3it)$$

$$(13)\times(12):
\f{e^{(-\f 12-3it-2\e,\e,\f 12+3it+\e)(c,0,-c)}}{
(-\f 12-3it-3\e)(-\f 12-3it)} R(\f 12+3it+3\e) R(-\f 12+3it+3\e) $$

$$(13)\times(13):
\f{e^{(-2\e,-1+\e,1+\e)(c,0,-c)}}{(1-3\e)(-2)}
R(\f 12+3it+3\e) R(-\f 12+3it+3\e)
R(\f 12-3it) R(-\f 12-3it) $$

$$(13)\times(321):
\f{e^{(-\f 12-3it-2\e,-\f 12+3it+\e ,1+\e)(c,0,-c)}}{
(-6it-3\e)(-\f 32+3it)} R(\f 12+3it+3\e) R(-\f 12+3it+3\e)
 R(-\f 12-3it) $$

$$(321)\times(12):
\f{e^{(-1+\e,\f 12-3it-2\e,\f 12+3it+\e)(c,0,-c)}}{
(-\f 32+3it+3\e)(-6it-3\e)} R(\f 12+3it+3\e)$$

$$(321)\times(13):
\f{e^{(-\f 12+3it+\e,-\f 12-3it-2\e,1+\e)(c,0,-c)}}{
(6it+3\e)(-\f 32-3it-3\e)}R(\f 12+3it+3\e) R(\f 12-3it)
R(-\f 12-3it)$$

$$(321)\times(321):
\f{e^{(-1+\e,-2\e,-1+\e)(c,0,-c)}}{(-1+3\e)(-1-3\e)} R(\f
12+3it+3\e) R(\f 12-3it).$$ The next two terms had limits taken in
$\d_1,\d_2\rightarrow 0$ and, up to constants from the residues,
are $$(12)\times e: \f{e^{(\e,\e,-2\e)(c,0,-c)}}{3\e}$$

$$(13)\times(123):
\f{e^{(-2\e,\e,\e)(c,0,-c)}}{-3\e} R(\f 12+3it+3\e) R(-\f 12+3it+3\e)
R(\f 12-3it) R(-\f 12-3it).$$

Note that
$$R(-\f 12+x) R(\f 12+x)=\f{Z(-\f 12+x)}{Z(\f 32+x)}
= \f{Z(\f 32-x)}{Z(\f 32+x)} $$
by the functional equation.
Accordingly, we may set $\e=0$ in the first 9 terms;
the last two involve a limit at $\e=0$ producing
derivatives.
%
%
%
%

\begin{prop}
\begin{equation}(i) ~~~ R(\f 12+it) \ll
1\label{pseveighti}\end{equation}

\begin{equation}(ii) ~~~~ \left|R(-\f{1}{2}+it) R(\f
12+it)\right|=1\label{pseveightii}\end{equation}

\begin{equation}(iii) ~~~~ \f{d}{d \e} \f{Z(\f 32+it+\e)}{Z(\f
32+it-\e)}\mid_{\e=0}=O(\log t).\label{pseveightiii}\end{equation}
\end{prop}
{\bf Proof}: (ii) follows from the facts that
$Z(\bar{s})=\overline{Z(s)}$ and $R(s)R(-s)=1$.

For (i), recall $$R(s)=\sqrt{\pi}\f{\G(\f{s}{2})}{\G(\f{s+1}{2})}
\f{\zeta(s)}{\zeta(s+1)}.$$  By Stirling's formula
$$\left|\f{\G(\f{1}{4}+\f{it}{2})}{\G(\f{3}{4}+\f{it}{2})}\right|\sim
\sqrt{\f{2}{t}}\mbox{ as } t\rightarrow \infty,$$ and by the
 convexity bound,
\begin{equation}\zeta(\f{1}{2}+it)=O_\e(t^{1/4+\e})\label{weylconvx}\end{equation}
 (better bounds can be
obtained -- see \cite{Titchmarsh}).  Also, taking the logarithm of
the Euler product of $\zeta(s)$ we find $$-\log
\zeta(\f{3}{2}+it)=\sum_{p\mbox{ prime}}\log(1-p^{-3/2-it}) =
O(1),$$ which proves (i).  Differentiating,
$$\f{\zeta'}{\zeta}(\f{3}{2}+it) =
-\sum_{p}\f{p^{-3/2-it}\log{p}}{1-p^{-3/2-it}}=O(1)$$ also.  So
$$\f{Z'}{Z}(\f{3}{2}+it)=-\f{1}{2}\log\pi +
\f{1}{2}\f{\G'}{\G}(\f{3}{4}+\f{it}{2}) +
 \f{\zeta'}{\zeta}(\f{3}{2}+it)
 = O(\log t),$$
proving (iii).  \bx

Summarizing these pointwise bounds,
\begin{prop}
\label{maxconstbd} For $T$ large, $$\int_{T-1}^{T+1}\int_{\F}|\L^C
E(P_1,g,(it,it,-2it),1)|^2dgdt =O_{\e}(T^\e).$$
\end{prop}

\subsection*{Assembling the Bounds}

{\bf Proof of (\ref{eisbound}):} Choose $c$ large enough so that
$\F_C$ is contained in $\{x \mid \hat{\tau}_{P_1}(H_0(x-C)),
\hat{\tau}_{P_2}(H_0(x-C))\le 0\}$.  Then the truncation does not
affect $\F_C$ and $$\int_{\F_C}|E(g,\l)|^2dg \le
\int_{\F}|\L^CE(g,\l)|^2dg$$ for each Eisenstein series $E(g,\l)$
on $\H$.  In the propositions we bounded local integrals of all of
the Eisenstein series as growing slower than any polynomial, with
the exception of the maximal parabolic Eisenstein series
(Remark~\ref{specsum}). For this we must sum over the
$SL_2({\Z})\backslash \U$ spectrum as well, so the contribution of
the left-hand side of (\ref{maxpbd}) near the point $\l$ is
bounded by $O_{\e}(\|\l\|^{1+\e})$.

 We note in comparison with \propref{gtandchar} that (\ref{paley})
 shows
$$ \hat{g}_{t\Sigma}(\l)=  \begin{cases}
    1+O_m(\mbox{dist}(\l,\partial (t\Sigma))^{-m}), &  \l \in  t\Sigma, \\
    O_m(\mbox{dist}(\l,\partial (t\Sigma))^{-m}), &  \l \in  t\Sigma
  \end{cases}, ~~~~~m\ge 0.
$$  Since we have shown that all the Eisenstein series grow
polynomially, we may switch the order of integration, ignore the
tail, and conclude (see (\ref{bsig}))
$$\int_{\F_C}Eis_{t\Sigma}(x,x)dx \ll
\int_{t\Sigma}(1+\|\l\|)^{1+\e} d\l.$$   Since $\b(\l)$ grows at
the rate of $\|\l\|^3$ in almost all directions, we have completed
the proof that $$\int_{\F_C} Eis_{g_{t\Sigma}}(x,x)dx =
o\(\int_{t\Sigma} \b(\l)d\l \).$$ \bx

\end{document}